\documentstyle[12pt]{article}

\newfont{\calmin}{cmff10}
\newcommand{\g}{\mbox{\calmin g}}
\def\GLh{$\mbox{GL}_h(n) \times \mbox{GL}_h(m)$}
\def\GLq{$\mbox{GL}_q(n) \times \mbox{GL}_q(m)$}
\def\cR{\mbox{$\cal{R}$}}
\def\cT{\mbox{$\cal{T}$}}
\def\cC{\mbox{$\cal{C}$}}
\def\cI{\mbox{$\cal{I}$}}
\def\bA{\mbox{\boldmath $A$}}
\def\bI{\mbox{\boldmath $I$}}
\def\bC{\mbox{\boldmath $C$}}
\def\bg{\mbox{\boldmath $g$}}
\def\tbA{\mbox{\hskip 4pt $\tilde{\mbox{\hskip -4pt \bA}}$}}
\def\tR{\tilde{R}}
\def\tcR{\mbox{\hskip 2pt $\tilde{\mbox{\hskip -2pt ${\cal R}$}}$}}
\def\tT{\tilde{T}}
\def\tcT{\mbox{\hskip 2pt $\tilde{\mbox{\hskip -2pt ${\cal T}$}}$}}
\def\tA{\mbox{\hskip 4pt $\tilde{\mbox{\hskip -4pt $A$}}$}}
\def\case#1#2{{\textstyle{#1\over #2}}}

\begin{document}

\title{\hfill{\normalsize ULB/229/CQ/98/4}\\
\vspace{1cm}
\boldmath Nonstandard GL$_h$($n$) quantum groups and contraction of
covariant $q$-bosonic algebras\footnote{Presented at the 7th Colloquium
``Quantum
Groups and Integrable Systems'', Prague, 18--20 June 1998}}
\author{C. Quesne\thanks{Directeur de recherches FNRS; E-mail:
cquesne@ulb.ac.be}\\
{\small \sl Physique Nucl\'eaire Th\'eorique et Physique Math\'ematique,
Universit\'e Libre de Bruxelles,}\\
{\small \sl Campus de la Plaine CP229, Boulevard du
Triomphe, B-1050 Brussels, Belgium}}
%
%
\maketitle

\begin{abstract}
\GLh-covariant $h$-bosonic algebras are built by contracting the \GLq-covariant
$q$-bosonic algebras considered by the present author some years ago. Their
defining relations are written in terms of the corresponding $R_h$-matrices.
Whenever $n=2$, and $m=1$ or~2, it is proved by using U$_h$(sl(2))
Clebsch-Gordan
coefficients that they can also be expressed in terms of coupled
commutators in a
way entirely similar to the classical case. Some U$_h$(sl(2)) rank-1/2
irreducible
tensor operators, recently contructed by Aizawa in terms of standard bosonic
operators, are shown to provide a realization of the $h$-bosonic algebra
corresponding to $n=2$ and $m=1$.
\end{abstract}

%
%
\section{Introduction}
It is well known that the Lie group GL(2) admits, up to isomorphism, only two
quantum group deformations with central determinant: the standard
deformation~GL$_q$(2), and the Jordanian deformation GL$_h$(2)~\cite{kuper}. The
quantum group~GL$_h$(2), or SL$_h$(2), and the dual quantum algebra of the
latter,
U$_h$(sl(2))~\cite{ohn}, have been the subject of many recent investigations,
among which one may quote the determination of the U$_h$(sl(2)) universal
$\cR$-matrix~\cite{ballesteros}.\par
%
%
Two useful tools have been devised for the Jordanian deformation study. One
of them is a contraction procedure that allows one to construct the latter
from the
standard deformation~\cite{agha}. In other words, GL$_h$(2) can be obtained from
GL$_q$(2) by a singular limit of a similarity transformation. Such a
technique has
been generalized by Alishahiha to higher-dimensional quantum
groups~\cite{ali}.\par
%
%
The other tool is a nonlinear invertible map between the generators of
U$_h$(sl(2))
and sl(2)~\cite{abdesselam}, yielding an explicit and simple method for
constructing the finite-dimensional irreducible representations~(irreps) of
U$_h$(sl(2)). In addition, it has provided an explicit formula for U$_h$(sl(2))
Clebsch-Gordan coefficients (CGC)~\cite{joris}, as well as bosonic or fermionic
realizations of irreducible tensor operators~(ITO) for
U$_h$(sl(2))~\cite{aizawa}.\par
%
%
The purpose of the present communication is to apply the contraction
procedure of
Ref.~\cite{agha}, as generalized by Alishahiha~\cite{ali}, to the \GLq-covariant
$q$-bosonic algebras constructed by the present author some years
ago~\cite{cq}, and recently rederived by Fiore by another
procedure~\cite{fiore98}.
As a result, we will obtain \GLh-covariant $h$-bosonic algebras. We will then
consider the cases where $n=2$, $m=1$, and $n=m=2$ in more detail, and establish
some relations with the works of Aizawa on ITO~\cite{aizawa}, and of Van der
Jeugt on CGC for U$_h$(sl(2))~\cite{joris}.\par
%
%
\section{\boldmath Contraction of GL$_q$($N$)}
The quantum group~GL$_q$($N$) is defined by the $RTT$-relations, $R' T'_1 T'_2 =
T'_2 T'_1 R'$, where $T' = \left(T'_{ij}\right) \in \mbox{GL}_q(N)$, $T'_1 = T'
\otimes I$, $T'_2 = I \otimes T'$, and
\begin{equation}
  R' = R'_q = q \sum_i e_{ii} \otimes e_{ii} + \sum_{i\ne j} e_{ii} \otimes
e_{jj}
  + \left(q - q^{-1}\right) \sum_{i<j} e_{ij} \otimes e_{ji},
\end{equation}
with $i$, $j$ running over 1, 2, \ldots,~$N$, and $e_{ij}$ denoting the $N
\times
N$~matrix with entry~1 in row~$i$ and column~$j$, and zeros everywhere else. An
equivalent form of the $RTT$-relations is obtained by replacing $R' = R'_{12}$
by~$R'^{-1}_{21}$. Throughout this communication, $q$-deformed objects
will be denoted by primed quantities, whereas unprimed ones will represent
$h$-deformed objects.\par
%
%
Let us consider the similarity transformation $R'' = \left(g^{-1} \otimes
g^{-1}\right) R' (g \otimes g)$, $T'' = g^{-1} T' g$, where $g$ is the $N
\times N$
matrix defined by $g = \sum_i e_{ii} + \eta e_{1N}$, in terms of some parameter
$\eta = h/(q-1)$~\cite{agha,ali}. The $RTT$-relations simply become $R''
T''_1 T''_2
= T''_2 T''_1 R''$.\par
%
%
Whenever $q$ goes to~1, although $\eta$ becomes singular, the latter have a
definite limit $R T_1 T_2 = T_2 T_1 R$, where $T = \lim_{q\to1} T''$, and
\begin{eqnarray}
  R & = & R_h = \lim_{q\to1} R'' \nonumber \\
  & = & \sum_{ij} e_{ii} \otimes e_{jj} + h \biggl[e_{11} \otimes e_{1N} -
e_{1N}
       \otimes e_{11} + e_{1N} \otimes e_{NN} - e_{NN} \otimes e_{1N}
\nonumber \\
  & & \mbox{} + 2 \sum_{i=2}^{N-1} (e_{1i} \otimes e_{iN} - e_{iN} \otimes
e_{1i})
       \biggr] + h^2 e_{1N} \otimes e_{1N}.
\end{eqnarray}
The resulting $R$-matrix is triangular, i.e., it is quasitriangular and
$R_{12}^{-1} =
R_{21}$, showing that the two equivalent forms of $RTT$-relations for
GL$_q$($N$)
have actually the same contraction limit. The matrix elements $T_{ij}$ generate
GL$_h$($N$).\par
%
%
\section{\boldmath \GLq-covariant $q$-bosonic algebras}
Let us consider two different copies of~GL$_q$($N$), corresponding to possibly
different dimensions $n$,~$m$, and let us denote quantities referring
to~GL$_q$($n$) by ordinary letters ($R'$, $T'$,~\ldots), and quantities referring
to~GL$_q$($m$) by script ones ($\cR'$, $\cT'$,~\ldots). The elements
$T'_{ij}$, $i$,
$j=1$, 2,~$\ldots n$, of GL$_q$($n$), and $\cT'_{st}$, $s$, $t=1$,
2,~$\ldots m$, of
GL$_q$($m$) are assumed to commute with one another.\par
%
%
In Ref.~\cite{cq}, $q$-bosonic creation and annihilation operators
$\bA'^+_{is}$, $\tbA'_{is}$, $i=1$,  2,~$\ldots n$, $s=1$, 2,~$\ldots m$,
that are
double ITO of rank $[1 \dot{0}]_n [1 \dot{0}]_m$, and $[\dot{0} -1]_n
[\dot{0}-1]_m$
with respect to $\mbox{U}_q(gl(n)) \times \mbox{U}_q(gl(m))$, respectively, were
constructed in terms of standard $q$-bosonic operators~\cite{biedenharn}
$a'^+_{is}$, $a'_{is}$, $i=1$, 2, \ldots,~$n$, $s=1$, 2, \ldots,~$m$,
acting in a tensor
product Fock space $F = \prod_{i=1}^n \prod_{s=1}^m F_{is}$. The annihilation
operators $\bA'_{is}$ contragredient to $\bA'^+_{is}$ were also considered. Both
sets of annihilation operators $\tbA'_{is}$ and $\bA'_{is}$, $i=1$, 2,
\ldots,~$n$,
$s=1$, 2, \ldots,~$m$, are related through the equation $\tbA' = \bA'
\bC'$, where
$\bC' = C' \cC'$, $C' = \sum_i (-1)^{n-i} q^{-(n-2i+1)/2} e_{ii'}$, and
$\cC' = \sum_s
(-1)^{m-s} q^{-(m-2s+1)/2} e_{ss'}$, with $i' = n-i+1$, $s' = m-s+1$.\par
%
%
The operators $\bA'^+_{is}$, $\bA'_{is}$, or $\bA'^+_{is}$, $\tbA'_{is}$,
generate with $\bI = I \cI$ a $\mbox{U}_q(\mbox{gl}(n)) \times
\mbox{U}_q(\mbox{gl}(m))$-module algebra or \GLq-comodule algebra, whose
$q$-commutation relations can be compactly written in coupled form by using
$\mbox{U}_q(\mbox{gl}(n)) \times \mbox{U}_q(\mbox{gl}(m))$ CGC. When rewritten
in componentwise form, such relations can be expressed in terms of the
GL$_q$($n$) and GL$_q$($m$) $R$-matrices as~\cite{cq}
\begin{eqnarray}
  R' \bA'^+_1 \bA'^+_2 & = & \bA'^+_2 \bA'^+_1 \cR', \qquad R' \bA'_2
\bA'_1 = \bA'_1
         \bA'_2 \cR', \nonumber \\
  \bA'_2 \bA'^+_1 & = & \bI_{21} + R'^{t_1} \cR'^{t_1} \bA'^+_1 \bA'_2,
         \label{eq:R'-first}
\end{eqnarray}
or
\begin{eqnarray}
  R' \bA'^+_1 \bA'^+_2 & = & \bA'^+_2 \bA'^+_1 \cR', \qquad R' \tbA'_1 \tbA'_2 =
         \tbA'_2 \tbA'_1 \cR', \nonumber \\
  \tbA'_2 \tbA'^+_1 & = & \bC'_{12} + q^2 \bA'^+_1 \tbA'_2 \tR'^{-1} \tcR'^{-1},
         \label{eq:R'-second}
\end{eqnarray}
where $t_1$ (resp.~$t_2$) denotes transposition in the first (resp.~second)
space
of the tensor product, $\tR'$ is defined by $\tR' = q C'_1
\left(R'^{-1}\right)^{t_1}
C'^{-1}_1 = q C'_2 \left(R'^{t_2}\right)^{-1} C'^{-1}_2$, and similar
relations hold
for~$\tcR'$. The transformations leaving Eqs.~(\ref{eq:R'-first})
and~(\ref{eq:R'-second}) invariant are $\varphi'\left(\bA'^+\right) =
\bA'^+ T' \cT'$,
$\varphi'(\bA') = T'^{-1} \cT'^{-1} \bA'$, and $\varphi'\left(\bA'^+\right)
= \bA'^+ T'
\cT'$, $\varphi'\bigl(\tbA'\bigr) = \tbA' \tT' \tcT'$, respectively. Here
$\tT'$ and
$\tcT'$ are defined by $\tT' = C'^{-1} \left( T'^{-1}\right)^t C'$, and
$\tcT' = \cC'^{-1}
\left(\cT'^{-1}\right)^t \cC'$.\par
%
%
There exists another independent set of \GLq-covariant $q$-bosonic operators,
which satisfy equations similar to Eq.~(\ref{eq:R'-first})
or~(\ref{eq:R'-second}),
but with $R'_{12} \to R'^{-1}_{21}$, $\cR'_{12} \to \cR'^{-1}_{21}$,
implying $q^{-1}
\tR'_{12} \to q \tR'^{-1}_{21}$, $q^{-1} \tcR'_{12} \to q \tcR'^{-1}_{21}$.\par
%
%
\section{\boldmath \GLh-covariant $h$-bosonic algebras}
Let us apply the contraction procedure of Sec.~2 to the \GLq-covariant
$q$-bosonic
algebras, given in two equivalent forms in Eqs.~(\ref{eq:R'-first})
and~(\ref{eq:R'-second}), respectively. Since we now have two copies of
GL$_q$($N$), we have to consider two transformation matrices $g = \sum_i
e_{ii} +
\eta e_{1n}$, and $\g = \sum_s e_{ss} + \eta e_{1m}$, acting on GL$_q$($n$) and
GL$_q$($m$), respectively.\par
%
%
Let us first consider Eq.~(\ref{eq:R'-first}), and introduce transformed
$q$-bosonic
operators defined by $\bA''^+ = \bA'^+ \bg$, $\bA'' = \bg^{-1} \bA'$, where
$\bg = g\:
\g$. By using the property $R'^t_{12} = R'_{21}$, and a similar one for
$\cR'$, it is
straightforward to show that Eq.~(\ref{eq:R'-first}) becomes
\begin{eqnarray}
  \bA''^+_1 \bA''^+_2 & = & \bA''^+_2 \bA''^+_1 R''^{-1}_{21} \cR''_{12}, \qquad
         \bA''_1 \bA''_2 = R''_{12} \cR''^{-1}_{21} \bA''_2 \bA''_1,
\nonumber \\
  \bA''_2 \bA''^+_1 & = & \bI_{21} + R''^{t_1} \cR''^{t_1} \bA''^+_1 \bA''_2.
\end{eqnarray}
Since $R$ and $\cR$ are triangular, in the $q\to1$ limit the $h$-bosonic
operators
$\bA^+_{is} = \lim_{q\to1} \bA''^+_{is}$, $\bA_{is} = \lim_{q\to1} \bA''_{is}$
satisfy the relations
\begin{eqnarray}
  \bA^+_1 \bA^+_2 & = & \bA^+_2 \bA^+_1 R \cR, \qquad \bA_1 \bA_2 = R \cR \bA_2
          \bA_1, \nonumber \\
  \bA_2 \bA^+_1 & = & \bI_{21} + R^{t_1} \cR^{t_1} \bA^+_1 \bA_2,
          \label{eq:R-first}
\end{eqnarray}
defining a \GLh-comodule algebra. The transformation
$\varphi\left(\bA^+\right) =
\bA^+ T \cT$, $\varphi(\bA) = T^{-1} \cT^{-1} \bA$, where $T_{ij} \in \mbox{\rm
GL}_h(n)$, $\cT_{st} \in \mbox{\rm GL}_h(m)$, leaves Eq.~(\ref{eq:R-first})
invariant.\par
%
%
Three properties of Eq.~(\ref{eq:R-first}) are worth noting: (1)~Had we started
instead from the second form of Eq.~(\ref{eq:R'-first}) corresponding to the
substitutions $R'_{12} \to R'^{-1}_{21}$, $\cR'_{12} \to \cR'^{-1}_{21}$,
we would
have obtained the same contraction limit~(\ref{eq:R-first}), owing to the
triangularity of~$R$ and~$\cR$. (2)~Contrary to what happens in the $q$-bosonic
case, $\bA_{is}$ can never be considered as the adjoint of $\bA^+_{is}$,
since no
*-structure is known on~GL$_h$($N$). (3)~For $m=1$, Eq.~(\ref{eq:R-first}) is
consistent with the general form of $\cal H$-covariant deformed bosonic algebras
for triangular~$\cal H$, obtained by Fiore~\cite{fiore97}.\par
%
%
Let us next consider Eq.~(\ref{eq:R'-second}), and define $\bA''^+ = \bA'^+
\bg$,
$\tbA'' = \tbA' \bg$, where $\bg$ is the same as before. Compatibility of the
$\tbA''$ and $\bA''$ definitions with $\tbA'' = \bA'' \bC''$, where $\bC''
= C'' \cC''$,
leads to $C'' = g^t C' g$, $\cC'' = \g^t \cC' \g$. A simple calculation
shows that for
$n>1$, a contraction limit of~$C''$ only exists for even $n$~values, and is
given by
$C = \lim_{q\to1} C'' = \sum_i (-1)^i e_{ii'} + (n-1) h e_{nn}$. Similar
results hold
for $\cC = \lim_{q\to1} \cC''$.\par
%
%
Restricting the range of $n$, $m$ values to $\{1, 2, 4, 6, \ldots\}$, we
obtain that
after transformation, Eq.~(\ref{eq:R'-second}) contracts into
\begin{eqnarray}
  \bA^+_1 \bA^+_2 & = & \bA^+_2 \bA^+_1 R \cR, \qquad \tbA_1 \tbA_2 =
       \tbA_2 \tbA_1 R \cR, \nonumber \\
  \tbA_2 \bA^+_1 & = & \bC_{12} + \bA^+_1 \tbA_2 \tR^{-1} \tcR^{-1},
       \label{eq:R-second}
\end{eqnarray}
where $\bC = C \cC$, $\tR = \lim_{q\to1} \left(g^{-1} \otimes g^{-1}\right)
\tR' (g
\otimes g) = C^{-1}_1 \left(R^{-1}\right)^{t_1} C_1 = C^{-1}_2
\left(R^{t_2}\right)^{-1} C_2$, and similarly for $\tcR$. For such
restricted $n$,
$m$ values, Eq.~(\ref{eq:R-second}) yields another form of the \GLh-covariant
$h$-bosonic algebra defined in Eq.~(\ref{eq:R-first}) for arbitrary $n$,
$m$ values.
The transformation leaving Eq.~(\ref{eq:R-second}) invariant is
$\varphi\left(\bA^+\right) = \bA^+ T \cT$, $\varphi\bigl(\tbA\bigr) = \tbA \tT
\tcT$, where $\tT = C^{-1} \left(T^{-1}\right)^t C$, $\tcT = \cC^{-1}
\left(\cT^{-1}\right)^t \cC$. However, for $n$ and/or $m \in \{3, 5, 7,
\ldots\}$,
the contraction procedure does not preserve the equivalence between
Eqs.~(\ref{eq:R'-first}) and~(\ref{eq:R'-second}), since only the former has a
limit.\par
%
%
\section{\boldmath GL$_h$(2) and $\mbox{GL}_h(2) \times
\mbox{GL}_h(2)$-covariant $h$-bosonic algebras}
{}For $n=2$, $m=1$, by making the substitutions
\begin{equation}
  R = \left(\begin{array}{cccc}
        1 & h & -h & h^2 \\[0.1cm]
        0 & 1 & 0  & h \\[0.1cm]
        0 & 0 & 1  & -h \\[0.1cm]
        0 & 0 & 0  & 1
        \end{array}\right), \qquad
  C = \left(\begin{array}{cc}
        0 & -1 \\[0.1cm]
        1 & h
        \end{array}\right), \qquad
  \cR = \cC = 1,   \label{eq:n=2,m=1}
\end{equation}
into Eqs.~(\ref{eq:R-first}) and~(\ref{eq:R-second}), we obtain that $A^+_1$,
$A^+_2$, $A_1$, $A_2$ satisfy the commutation relations
\begin{eqnarray}
  \left[A^+_1, A^+_2\right] & = & h \left(A^+_1\right)^2, \qquad \left[A_1,
        A_2\right] = h A^2_2, \nonumber \\
  \left[A_2, A^+_1\right] & = & 0, \qquad \left[A_1, A^+_2\right] = h
\left(- A^+_1
        A_1 - A^+_2 A_2 + h A^+_1 A_2\right), \nonumber \\
  \left[A_1, A^+_1\right] & = & \left[A_2, A^+_2\right] = I + h A^+_1 A_2,
        \label{eq:commut-1}
\end{eqnarray}
while $A^+_1$, $A^+_2$, $\tA_1$, $\tA_2$ fulfil
\begin{eqnarray}
  \left[A^+_1, A^+_2\right] & = & h \left(A^+_1\right)^2, \qquad \bigl[\tA_1,
        \tA_2\bigr] = h \tA^2_1, \nonumber \\
  \bigl[\tA_1, A^+_1\bigr] & = & 0, \qquad \bigl[\tA_2, A^+_2\bigr] = h
\bigl(I -
        A^+_1 \tA_2 + A^+_2 \tA_1 + h A^+_1 \tA_1\bigr), \nonumber \\
  \bigl[\tA_1, A^+_2\bigr] & = & - \bigl[\tA_2, A^+_1\bigr] = I + h A^+_1 \tA_1.
        \label{eq:commut-2}
\end{eqnarray}
\par
%
%
Both sets of operators $\left(A^+_1, A^+_2\right)$ and $\bigl(\tA_1,
\tA_2\bigr)$
may be considered as the components $m = 1/2$ and $m=-1/2$ of ITO of rank~1/2,
or spinors, with respect to the quantum algebra~U$_h$(sl(2)). By considering the
adjoint action of the U$_h$(sl(2)) generators on such spinors,
Aizawa~\cite{aizawa} recently realized them in terms of standard bosonic
operators $a^+_1$, $a^+_2$, $a_1$, $a_2$,
\begin{eqnarray}
  A^+_1 & = & \left(1 - \case{h}{2} J_+\right)^{-1} a^+_1, \qquad A^+_2 =
\left(1 -
         \case{h}{2} J_+\right) a^+_2 + \case{h}{2} \left(A^+_1 - 2 a^+_1
J_0\right),
         \nonumber \\
  \tA_1 & = & \left(1 - \case{h}{2} J_+\right)^{-1} a_2, \qquad \tA_2 = -
\left(1 -
         \case{h}{2} J_+\right) a_1 + \case{h}{2} \left(\tA_1 - 2 a_2
J_0\right),
         \label{eq:realization}
\end{eqnarray}
where $J_+ = a^+_1 a_2$, and $J_0 = \left(a^+_1 a_1 - a^+_2 a_2\right)/2$
are sl(2)
generators. As can be easily checked, the operators~(\ref{eq:realization})
satisfy
Eq.~(\ref{eq:commut-2}), as it should be.\par
%
%
Equation~(\ref{eq:commut-2}) can be recast into an alternative form by using
coupled commutators
\begin{equation}
  \left[U^{j_1}, V^{j_2}\right]^j_m \equiv \left[U^{j_1} \times
V^{j_2}\right]^j_m
  - (-1)^{\epsilon} \left[V^{j_2} \times U^{j_1}\right]^j_m,
\end{equation}
where $U^{j_1}$ and $V^{j_2}$ denote two ITO of rank $j_1$ and $j_2$ with
respect to U$_h$(sl(2)), respectively, $\epsilon = j_1 + j_2 - j$,
\begin{equation}
  \left[U^{j_1} \times V^{j_2}\right]^j_m \equiv \sum_{m_1m_2} \langle j_1
  m_1, j_2 m_2 | j m \rangle_h\, U^{j_1}_{m_1} V^{j_2}_{m_2},
\end{equation}
and $\langle\, ,\, | \,\rangle_h$ denotes a U$_h$(sl(2)) CGC, as determined in
Ref.~\cite{joris}. The results read
\begin{equation}
  \bigl[A^+, A^+\bigr]^0_0 = \bigl[\tA, \tA\bigr]^0_0 = \bigl[\tA,
A^+\bigr]^1_m = 0,
  \qquad \bigl[\tA, A^+\bigr]^0_0 = \sqrt{2}\, I.  \label{eq:coupled-1}
\end{equation}
\par
%
%
{}For $n=m=2$, $\cR$ and $\cC$ take the same form as $R$ and $C$ in
Eq.~(\ref{eq:n=2,m=1}). Relations similar to those in Eqs.~(\ref{eq:commut-1})
and~(\ref{eq:commut-2}) can be easily written. The operators $\bA^+_{is}$,
$\tbA_{is}$, $i$, $s=1$,~2, may now be considered as the components of double
spinors with respect to $\mbox{U}_h(\mbox{sl}(2)) \times
\mbox{U}_h(\mbox{sl}(2))$, and they satisfy the coupled commutation relations
\begin{eqnarray}
  \bigl[\bA^+, \bA^+\bigr]^{1,0}_{m,0} & = & \bigl[\bA^+,
\bA^+\bigr]^{0,1}_{0,m'} =
          \bigl[\tbA, \tbA\bigr]^{1,0}_{m,0} = \bigl[\tbA,
\tbA\bigr]^{0,1}_{0,m'} = 0,
          \nonumber \\
  \bigl[\tbA, \bA^+\bigr]^{j,j'}_{m,m'} & = & 2 \delta_{j,0} \delta_{j',0}
\delta_{m,0}
          \delta_{m',0} \bI,  \label{eq:coupled-2}
\end{eqnarray}
where in the definition of coupled commutators there now appear two
$\epsilon$~phases, and two U$_h$(sl(2)) CGC.\par
%
%
It is remarkable that both Eqs.~(\ref{eq:coupled-1})
and~(\ref{eq:coupled-2}) are
formally identical with those for sl(2) and $\mbox{sl(2)} \times \mbox{sl(2)}$,
respectively. Contrary to what happens in the $q$-bosonic case where the
commutators are $q$-deformed, here all the dependence upon the deforming
parameter~$h$ is contained in the CGC.\par
%
%
\section{Conclusion}
In this communication, we showed that \GLh-covariant $h$-bosonic
algebras can be obtained by contracting \GLq-covariant $q$-bosonic ones. Some
extensions of the present work to $h$-fermionic and multiparametric algebras
are under current investigation.\par
%
%

\end{document}